\documentclass[12pt]{article}
\usepackage{latexsym}
\usepackage{amsfonts}
\usepackage{amsmath}
\usepackage{amssymb}
\usepackage{graphicx}
\usepackage{latexsym}
\usepackage{amsfonts}
\usepackage{graphicx}
\usepackage{psfrag}

\input{tcilatex}
\begin{document}

\qquad 

\thispagestyle{empty}

\begin{center}
{\Large \textbf{Some generalisations of the inequalities for positive linear
maps }}

\vskip0.2inR. Sharma, P. Devi and R. kumari

Department of Mathematics \& Statistics\\[0pt]
Himachal Pradesh University\\[0pt]
Shimla -5,\\[0pt]
India - 171005\\[0pt]
email: rajesh.sharma.hpn@nic.in
\end{center}

\vskip1.5in \noindent \textbf{Abstract. \ }We obtain generalisations of some
inequalities for positive unital linear maps on matrix algebra. This also
provides several positive semidefinite matrices, and we get some old and new
inequalities involving the eigenvalues of a Hermitian matrix.



\bigskip

\bigskip

\bigskip

\bigskip

\bigskip

\bigskip

\bigskip

\bigskip

\bigskip

\section{\protect\bigskip Introduction}

\setcounter{equation}{0}Let $\mathbb{M}\left( n\right) $ be the $C^{\ast }$-
algebra of all $n\times n$ complex matrices. Let $\Phi :\mathbb{M}\left(
n\right) \rightarrow \mathbb{M}\left( k\right) $ be a positive unital linear
map $\left[ 3\right] $. Kadison's inequality $\left[ 11\right] $ says that
for any Hermitian element $A$ of $\mathbb{M}\left( n\right) ,$ we have%
\begin{equation}
\Phi (A^{2})\geq \Phi (A)^{2}  \tag{1.1}
\end{equation}%
or equivalently%
\begin{equation}
\left[ 
\begin{array}{cc}
I & \Phi (A) \\ 
\Phi (A) & \Phi (A^{2})%
\end{array}%
\right] \geq O.  \tag{1.2}
\end{equation}%
For more details, generalisations and extensions of this inequality, see
Davis $\left[ 9\right] $ and Choi $\left[ 7,8\right] .$\vskip0.0inA
complementary inequality due to Bhatia and Davis $\left[ 2\right] $ says
that if the spectrum of a Hermitian matrix $A$ is contained in the interval $%
\left[ m,M\right] ,$ then%
\begin{equation}
\Phi (A^{2})-\Phi (A)^{2}\leq \left( \frac{M-m}{2}\right) ^{2}I.  \tag{1.3}
\end{equation}%
They also proved that 
\begin{equation}
\Phi (A^{2})-\Phi (A)^{2}\leq \left( \Phi (A)-mI\right) \left( MI-\Phi
(A)\right) .  \tag{1.4}
\end{equation}%
The inequality $\left( 1.4\right) $ provides a refinement of $\left(
1.3\right) .$ For more details and applications of these inequalities, see $%
\left[ 4-6\right] .$ \vskip0.0inThe Kadison inequality (1.1) is a
noncommutative version of the classical inequality%
\begin{equation}
\mathbb{E}\left( X^{2}\right) \geq \mathbb{E}\left( X\right) ^{2},  \tag{1.5}
\end{equation}%
where $\mathbb{E}\left( X\right) =\dint xf\left( x\right) dx$, $f\left(
x\right) \geq 0$ and $\dint f\left( x\right) dx=1$. Kadison $\left[ 11\right]
$ remarks that the standard proof of the corresponding inequalities for
scalars do not apply to give simple proofs for linear maps. In case of
Kadison's inequality $\left( 1.1\right) $, if $\Phi (A)$ and $\Phi (A^{2})$
commute one can reduce the problem to the real valued case, and the results
follow from these considerations. \vskip0.0inThe inequality $\left(
1.5\right) $ is subsumed in a more general Jensen's inequality that says
that if $f$ is a convex function on $\left( a,b\right) $ then%
\begin{equation}
f\left( \mathbb{E}(X)\right) \leq \mathbb{E}\left( f(X)\right) .  \tag{1.6}
\end{equation}%
One generalisation of the Kadison inequality (1.1) is a noncommutative
analogue of (1.6) that says that if $A$ is a Hermitian matrix whose spectrum
is contained in $\left( a,b\right) $ and if $f$ is matrix convex function $%
\left[ 3\right] $ on $\left( a,b\right) $ and $\Phi $ is a positive unital
linear map, then%
\begin{equation}
f\left( \Phi (A)\right) \leq \Phi \left( f(A)\right) ,  \tag{1.7}
\end{equation}%
see Davis $\left[ 9\right] $ and Choi $\left[ 7\right] .$ Bhatia and Sharma $%
[4]$ have shown that for $2\times 2$ matrices the inequality $\left(
1.7\right) $ holds true for all ordinary convex functions on $\left(
a,b\right) .$ It is well known that the same is true for every positive
unital linear functional $\varphi :\mathbb{M}(n)\rightarrow \mathbb{C}$.

\vskip0.0inThe inequality (1.5) is also subsumed in one more general
inequality that says that $\left[ \mathbb{E}\left( X^{i+j-2}\right) \right] $
is a real symmetric positive semidefinite matrix, see [10]. Here, we study
such possible extensions of the noncommutative inequalities $\left(
1.2\right) $ and $\left( 1.4\right) .$ The proof of the inequality $\left(
1.2\right) $ involves the Spectral theorem and the properties of tensor
product, see $\left[ 3\right] .$ By the Spectral theorem, 
\begin{equation*}
A=\dsum\limits_{k=1}^{n}\lambda _{k}P_{k},
\end{equation*}%
where $\lambda _{k}$ are the eigenvalues of $A$ and $P_{k}$ the
corresponding projections with $\dsum\limits_{k=1}^{n}P_{k}=I.$ Then 
\begin{equation}
A^{r}=\dsum\limits_{k=1}^{n}\lambda _{k}^{r}P_{k}\text{ , }\Phi \left(
A^{r}\right) =\dsum\limits_{k=1}^{n}\lambda _{k}^{r}\Phi \left( P_{k}\right) 
\text{ and }\dsum\limits_{k=1}^{n}\Phi \left( P_{k}\right) =I.  \tag{1.8}
\end{equation}%
We augment the technique of the proof $\left[ 3\right] $ of the inequality
(1.2) and prove one more generalisation of the Kadison inequality that says
that the matrix%
\begin{equation*}
\left[ 
\begin{array}{cccc}
I & \Phi (A) & \cdots & \Phi (A^{r}) \\ 
\Phi (A) & \Phi (A^{2}) & \cdots & \Phi (A^{r+1}) \\ 
\vdots & \vdots & \ddots & \vdots \\ 
\Phi (A^{r}) & \Phi (A^{r+1}) & \cdots & \Phi (A^{2r})%
\end{array}%
\right]
\end{equation*}%
is positive semidefinite, for all $r=0,1,2...\left( \text{see Corollary}%
\left( 2.1\right) ,\text{below}\right) .$ Our main results (Theorem $2.1-2.5$%
) provide some generalisations of the Bhatia-Davis inequality $\left(
1.4\right) .$ It is shown that some recent bounds for the central moments
are the special cases of our results. Lower and upper bounds respectively
for the largest and smallest eigenvalues of a Hermitian matrix are derived,
(Theorem 3.3). Inequalities analogous to Kadison's inequality (1.2) are
obtained (Theorem 3.2-3.3).

\section{Main results}

\setcounter{equation}{0}\textbf{Lemma 2.1. }For $x\geq y,$ the $(r+1)\times
(r+1)$ matrix $\left[ x^{i+j-1}-yx^{i+j-2}\right] $ is positive semidefinite
(psd), where $1\leq i,j\leq r+1,$ $r=0,1,2,...$.

\vskip0.1in\noindent \textbf{Proof. }The matrix $\left[ x^{i+j-2}\right] $
is psd for all $x\in \mathbb{R}$. Its rank is one and trace is non-negative.
Likewise, the matrix $\left[ x-y\right] $ is psd for $x\geq y$. The Schur
product

\begin{equation*}
\left[ x^{i+j-2}\right] \circ \left[ x-y\right] =\left[ x^{i+j-1}-yx^{i+j-2}%
\right] \text{ \ }
\end{equation*}%
is also psd. \ $\blacksquare $

\noindent \vskip0.1in\textbf{Theorem 2.1}.\textbf{\ }Let $\Phi :\mathbb{M}%
\left( n\right) \rightarrow \mathbb{M}\left( k\right) $ be a positive unital
linear map. Let $A$ be any Hermitian element of $\mathbb{M}\left( n\right) $
whose spectrum is contained in $\left[ m,M\right] .$ Then 
\begin{equation}
\left[ \Phi \left( A^{i+j-1}\right) -m\Phi \left( A^{i+j-2}\right) \right]
_{r+1\times r+1}\geq O  \tag{2.1}
\end{equation}%
and%
\begin{equation}
\left[ M\Phi \left( A^{i+j-2}\right) -\Phi \left( A^{i+j-1}\right) \right]
_{r+1\times r+1}\geq O,  \tag{2.2}
\end{equation}%
where $1\leq i,j\leq r+1,$ $r=0,1,2,..$.\newline
\noindent \textbf{Proof. }Using $\left( 1.8\right) ,$ we have\textbf{\ }%
\begin{equation}
\left[ \Phi \left( A^{i+j-1}\right) -m\Phi \left( A^{i+j-2}\right) \right]
=\sum_{k=1}^{n}\text{ }\left[ \left( \lambda _{k}^{i+j-1}-m\lambda
_{k}^{i+j-2}\right) \Phi \left( P_{k}\right) \right] .  \notag
\end{equation}%
Also,%
\begin{equation*}
\sum_{k=1}^{n}\text{ }\left[ \left( \lambda _{k}^{i+j-1}-m\lambda
_{k}^{i+j-2}\right) \Phi \left( P_{k}\right) \right] =\sum_{k=1}^{n}\text{ }%
\left[ \lambda _{k}^{i+j-1}-m\lambda _{k}^{i+j-2}\right] \otimes \Phi \left(
P_{k}\right)
\end{equation*}%
where $\otimes $ denotes the tensor product of matrices. So,%
\begin{equation}
\left[ \Phi \left( A^{i+j-1}\right) -m\Phi \left( A^{i+j-2}\right) \right]
=\sum_{k=1}^{n}\text{ }\left[ \lambda _{k}^{i+j-1}-m\lambda _{k}^{i+j-2}%
\right] \otimes \Phi \left( P_{k}\right) .  \tag{2.3}
\end{equation}%
Since $P_{k}$ is psd, $\Phi \left( P_{k}\right) $ are psd. It follows from
Lemma $2.1$ that the matrix $\left[ \lambda _{k}^{i+j-1}-m\lambda
_{k}^{i+j-2}\right] $ is psd. The tensor product of two psd matrices is psd.
Each summand in (2.3) is psd, and so is the sum. So, (2.1) holds true.
Likewise, (2.2) follows from the fact that the matrix $\left[ M\lambda
_{k}^{i+j-2}-\lambda _{k}^{i+j-1}\right] $ is psd. $\blacksquare \newline
$\textbf{Corollary 2.1.} Under the conditions of Theorem 2.1, 
\begin{equation}
\left[ \Phi \left( A^{i+j-2}\right) \right] _{r+1\times r+1}\geq O. 
\tag{2.4}
\end{equation}%
For $A\geq O,$ we also have%
\begin{equation}
\left[ \Phi \left( A^{i+j-1}\right) \right] _{r+1\times r+1}\geq O. 
\tag{2.5}
\end{equation}%
\textbf{Proof. }Adding $(2.1)$ and $(2.2)$, we immediately get $(2.4)$. The
inequality $(2.5)$ is a special case of $\left( 2.1\right) $; $m=0$. $%
\blacksquare $\newline

\vskip0.1inFor $\lambda _{k}\geq m>0,$ $k=1,2,...,n,$ the matrix $\left[ 1-%
\frac{m}{\lambda _{k}}\right] $ is psd and so is the matrix

\begin{equation}
\text{\ }\left[ \lambda _{k}^{i+j-2}\right] \circ \left[ 1-\frac{m}{\lambda
_{k}}\right] =\left[ \lambda _{k}^{i+j-2}-m\lambda _{k}^{i+j-3}\right] \text{
}.  \tag{2.6}
\end{equation}%
Using arguments similar to those used in the proof of Theorem 2.1, we have
the following theorem.

\noindent\ \textbf{\noindent Theorem 2}.\textbf{2.} Let $A,$ $\Phi ,m$ and $%
M $ be as in Theorem 2.1. If $A$ is positive definite, then%
\begin{equation}
\left[ \Phi \left( A^{i+j-2}\right) -m\Phi \left( A^{i+j-3}\right) \right]
_{r+1\times r+1}\geq O  \tag{2.7}
\end{equation}%
and%
\begin{equation}
\left[ M\Phi \left( A^{i+j-3}\right) -\Phi \left( A^{i+j-2}\right) \right]
_{r+1\times r+1}\geq O.\text{ \ }\blacksquare  \tag{2.8}
\end{equation}%
\textbf{Lemma 2.2. }For $x\geq y\geq z,$ the $(r+1)\times (r+1)$ matrix $%
\left[ y^{i+j-2}\left( x-y\right) \left( y-z\right) \right] $ is psd, $1\leq
i,j\leq r+1,$ $r=0,1,2,...$.

\vskip0.1in \noindent \textbf{Proof. }It is enough to note that the Schur
product 
\begin{equation}
\left[ y^{i+j-2}\left( x-y\right) \left( y-z\right) \right] =\left[ y^{i+j-2}%
\right] \circ \left[ \left( x-y\right) \left( y-z\right) \right]  \tag{2.9}
\end{equation}%
is psd. \ $\blacksquare $ \textbf{\noindent }\newline
\textbf{Theorem 2.3. }Under the condition of Theorem 2.1, we have%
\begin{equation}
\left[ \Phi (A^{i+j-2}(A-mI)(MI-A))\right] _{r+1\times r+1}\geq O. 
\tag{2.10}
\end{equation}%
\textbf{\noindent Proof}$.$\textbf{\ }Using the Lemma 2.2, we see that the
matrix%
\begin{equation*}
\text{\ }\left[ \lambda _{k}^{i+j-2}\right] \circ \left[ \left( \lambda
_{k}-m\right) \left( M-\lambda _{k}\right) \right] =\left[ \lambda
_{k}^{i+j-2}\left( \lambda _{k}-m\right) \left( M-\lambda _{k}\right) \right]
\end{equation*}%
is psd. The proof now follows on using arguments similar to those used in
the proof of Theorem 2.1. $\blacksquare $

\vskip0.1in\textbf{Theorem 2.4.} Let $\Phi :\mathbb{M}\left( n\right)
\rightarrow \mathbb{M}\left( k\right) $ be a positive unital linear map. Let 
$A$ be any Hermitian element of $\mathbb{M}\left( n\right) $ with distinct
eigenvalues $\lambda _{1}<\lambda _{2}<...<\lambda _{k}.$ Then 
\begin{equation}
\left[ \Phi (A^{i+j-2}(A-\lambda _{j-1}I)(A-\lambda _{j}I))\right]
_{r+1\times r+1}\geq O.  \tag{2.11}
\end{equation}%
\textbf{\ Proof}. Since all the $\lambda _{k}$ $(k=1,2,...,n)$ lies outsides 
$\left( \lambda _{j-1},\lambda _{j}\right) ,$ $j=2,3,...k$, we have $\left(
\lambda _{k}-\lambda _{j-1}\right) \left( \lambda _{k}-\lambda _{j}\right)
\geq 0.$ Therefore, the Schur product 
\begin{equation*}
\text{\ }\left[ \lambda _{k}^{i+j-2}\right] \circ \left[ \left( \lambda
_{k}-\lambda _{j-1}\right) \left( \lambda _{k}-\lambda _{j}\right) \right] =%
\left[ \lambda _{k}^{i+j-2}\left( \lambda _{k}-\lambda _{j-1}\right) \left(
\lambda _{k}-\lambda _{j}\right) \right] 
\end{equation*}%
is psd. The proof now follows on using arguments similar to those used in
the proof of Theorem 2.1. $\blacksquare $\vskip0.1in\textbf{Theorem 2.5.}
Let $A,$ $\Phi ,m$ and $M$ be as in Theorem 2.1. If $A$ is positive
definite, then%
\begin{equation}
\left[ \Phi \left( A^{i+j-3}(A-mI)(MI-A)\right) \right] _{r+1\times r+1}\geq
O.  \tag{2.12}
\end{equation}%
\textbf{\noindent Proof. }The arguments are similar to the proof of the
above theorem. Note that for $x\geq y\geq z$ and $y>0,$ we have $%
x+z-xzy^{-1}-y\geq 0.$ So, the Schur product 
\begin{equation*}
\text{\ }\left[ \lambda _{k}^{i+j-2}\right] \circ \left[ \left(
m+M-mM\lambda _{k}^{-1}-\lambda _{k}\right) \right] =\left[ \lambda
_{k}^{i+j-3}\left( \lambda _{k}-m\right) \left( M-\lambda _{k}\right) \right]
\end{equation*}%
is psd. $\blacksquare $

\section{Special Cases}

If $A$ and $B$ are positive definite matrices then the block matrix $\left[ 
\begin{array}{cc}
A & X \\ 
X^{\ast } & B%
\end{array}%
\right] $ is psd if and only if $A\geq XB^{-1}X^{\ast },$ see $\left[ 3%
\right] .$ Using this result we can discuss various special cases of
inequalities derived above. We demonstrate some of these cases here .

We find, on adding $\left( 2.7\right) $ and $\left( 2.8\right) $, that if $%
A>0$ then%
\begin{equation*}
\left[ 
\begin{array}{cccc}
\Phi (A^{-1}) & I & \cdots & \Phi (A^{r-1}) \\ 
I & \Phi (A) & \cdots & \Phi (A^{r}) \\ 
\vdots & \vdots & \ddots & \vdots \\ 
\Phi (A^{r-1}) & \Phi (A^{r}) & \cdots & \Phi (A^{2r-1})%
\end{array}%
\right]
\end{equation*}%
is psd. So,%
\begin{equation*}
\left[ 
\begin{array}{cc}
\Phi (A^{-1}) & I \\ 
I & \Phi (A)%
\end{array}%
\right] \geq O.
\end{equation*}%
This gives Choi $\left[ 7\right] $ inequality%
\begin{equation*}
\Phi (A^{-1})\geq \Phi (A)^{-1},
\end{equation*}%
see $\left[ 3\right] .$\vskip0.1inFor $r=1,$ Theorem 2.1 says that%
\begin{equation*}
\left[ 
\begin{array}{cc}
\Phi \left( A\right) -mI & \Phi \left( A^{2}\right) -m\Phi \left( A\right)
\\ 
\Phi \left( A^{2}\right) -m\Phi \left( A\right) & \Phi \left( A^{3}\right)
-m\Phi \left( A^{2}\right)%
\end{array}%
\right] \geq O
\end{equation*}%
and%
\begin{equation*}
\left[ 
\begin{array}{cc}
MI-\Phi \left( A\right) & M\Phi \left( A\right) -\Phi \left( A^{2}\right) \\ 
M\Phi \left( A\right) -\Phi \left( A^{2}\right) & M\Phi \left( A^{2}\right)
-\Phi \left( A^{3}\right)%
\end{array}%
\right] \geq O.
\end{equation*}%
For $\Phi (A)>mI,$ we therefore have%
\begin{equation}
\Phi \left( A^{3}\right) \geq m\Phi \left( A^{2}\right) +\left( \Phi \left(
A^{2}\right) -m\Phi \left( A\right) \right) \left( \Phi \left( A\right)
-mI\right) ^{-1}\left( \Phi \left( A^{2}\right) -m\Phi \left( A\right)
\right)  \tag{3.1}
\end{equation}%
and for $\Phi (A)<MI,$ we have%
\begin{equation}
\Phi \left( A^{3}\right) \leq M\Phi \left( A^{2}\right) -\left( M\Phi \left(
A\right) -\Phi \left( A^{2}\right) \right) \left( MI-\Phi \left( A\right)
\right) ^{-1}\left( M\Phi \left( A\right) -\Phi \left( A^{2}\right) \right) .
\tag{3.2}
\end{equation}%
For the corresponding commutative cases of the inequalities $\left(
3.1\right) $ and $\left( 3.2\right) ,$%
\begin{equation*}
\mathbb{E}\left( X^{3}\right) \geq m\mathbb{E}\left( X^{2}\right) +\frac{%
\left( \mathbb{E}\left( X^{2}\right) -m\mathbb{E}\left( X\right) \right) ^{2}%
}{\mathbb{E}\left( X\right) -m}
\end{equation*}%
and%
\begin{equation*}
\mathbb{E}\left( X^{3}\right) \leq M\mathbb{E}\left( X^{2}\right) -\frac{%
\left( M\mathbb{E}\left( X\right) -\mathbb{E}\left( X^{2}\right) \right) ^{2}%
}{M-\mathbb{E}\left( X\right) },
\end{equation*}%
see $\left[ 14\right] .$ The bounds for $\mathbb{E}\left( X^{2}\right) $ in
terms of $\mathbb{E}\left( X^{-1}\right) $ and $\mathbb{E}\left( X\right) $
are derived in $\left[ 13\right] $ on using derivatives. These inequalities
follow from our Theorem 2.2 for $r=1.$ Likewise, the bounds for $\mathbb{E}%
\left( X^{4}\right) $ in $\left[ 15\right] $ follow from our Theorem 2.3 for 
$r=1.$ The inequalities for the central moments $\left[ 13-15\right] $ also
follow from our more general results. Let $\varphi :\mathbb{M}(n)\rightarrow 
\mathbb{C}$ be a positive unital linear functional. Let $B=A-\varphi \left(
A\right) I$ , $a=m-\varphi \left( A\right) $ and $b=M-\varphi \left(
A\right) .$ It follows from Theorem 2.1 that%
\begin{equation}
\left[ \varphi \left( B^{i+j-1}\right) -a\varphi \left( B^{i+j-2}\right) %
\right] _{r+1\times r+1}\geq O  \tag{3.3}
\end{equation}%
and%
\begin{equation}
\left[ b\varphi \left( B^{i+j-2}\right) -\varphi \left( B^{i+j-1}\right) %
\right] _{r+1\times r+1}\geq O.  \tag{3.4}
\end{equation}%
From inequalities $\left( 3.3\right) $ and $\left( 3.4\right) $ , we also
have%
\begin{equation*}
\left[ \varphi \left( B^{i+j-2}\right) \right] _{r+1\times r+1}\geq O.
\end{equation*}%
Likewise, we can discuss the corresponding inequalities for functionals
related to Theorem 2.2-2.5.\vskip0.1inA special case of Corollary 2.1 says
that for every Hermitian matrix $A,$%
\begin{equation*}
\left[ 
\begin{array}{cc}
\Phi \left( A^{2}\right) & \Phi \left( A^{3}\right) \\ 
\Phi \left( A^{3}\right) & \Phi \left( A^{4}\right)%
\end{array}%
\right] \geq O.
\end{equation*}%
We prove a refinement of this inequality for positive definite matrices in
the following theorem.\vskip0.1in\textbf{Theorem 3.1. }Let $A,$ $\Phi $ and $%
m$ be as in Theorem 2.1. If $A\geq m>0,$ then%
\begin{equation}
\left[ 
\begin{array}{cc}
\Phi \left( A^{2}\right) & \Phi \left( A^{3}\right) \\ 
\Phi \left( A^{3}\right) & \Phi \left( A^{4}\right)%
\end{array}%
\right] \geq 2m\left[ 
\begin{array}{cc}
\Phi \left( A\right) & \Phi \left( A^{2}\right) \\ 
\Phi \left( A^{2}\right) & \Phi \left( A^{3}\right)%
\end{array}%
\right] -m^{2}\left[ 
\begin{array}{cc}
I & \Phi \left( A\right) \\ 
\Phi \left( A\right) & \Phi \left( A^{2}\right)%
\end{array}%
\right] \geq O.  \tag{3.5}
\end{equation}%
\textbf{Proof. }The second inequality $\left( 3.5\right) $ follows from the
fact that for $0<m\leq \lambda _{j}$, the Schur product%
\begin{equation*}
\left[ 
\begin{array}{cc}
m(2\lambda _{j}-m) & m\lambda _{j}(2\lambda _{j}-m) \\ 
m\lambda _{j}(2\lambda _{j}-m) & m\lambda _{j}^{2}(2\lambda _{j}-m)%
\end{array}%
\right] =\left[ 
\begin{array}{cc}
(2\lambda _{j}-m) & (2\lambda _{j}-m) \\ 
(2\lambda _{j}-m) & (2\lambda _{j}-m)%
\end{array}%
\right] \circ \left[ 
\begin{array}{cc}
m & m\lambda _{j} \\ 
m\lambda _{j} & m\lambda _{j}^{2}%
\end{array}%
\right]
\end{equation*}%
is psd.\vskip0.1inLikewise, the first inequality (3.5) follows from the fact
that the matrix%
\begin{equation*}
\left[ 
\begin{array}{cc}
(\lambda _{j}-m)^{2} & \lambda _{j}(\lambda _{j}-m)^{2} \\ 
\lambda _{j}(\lambda _{j}-m)^{2} & \lambda _{j}^{2}(\lambda _{j}-m)^{2}%
\end{array}%
\right] =\left[ 
\begin{array}{cc}
1 & \lambda _{j} \\ 
\lambda _{j} & \lambda _{j}^{2}%
\end{array}%
\right] \circ \left[ 
\begin{array}{cc}
(\lambda _{j}-m)^{2} & (\lambda _{j}-m)^{2} \\ 
(\lambda _{j}-m)^{2} & (\lambda _{j}-m)^{2}%
\end{array}%
\right]
\end{equation*}%
is psd. $\blacksquare $\vskip0.1inThe matrix%
\begin{equation*}
\left[ 
\begin{array}{cc}
\Phi \left( A^{2}\right) & \Phi \left( A\right) \\ 
\Phi \left( A\right) & \Phi \left( A\right)%
\end{array}%
\right]
\end{equation*}%
is not always psd. It is here interesting to note the following theorem.

\textbf{Theorem 3.2.} Let $\Phi $ be as in Theorem 2.1. For $A>0,$ we have 
\begin{equation*}
\left[ 
\begin{array}{cc}
\Phi \left( A^{2}\right) & \Phi \left( A\right) \\ 
\Phi \left( A\right) & \Phi \left( A-\log A\right)%
\end{array}%
\right] \geq O.
\end{equation*}%
\textbf{Proof}. We have 
\begin{equation*}
\left[ 
\begin{array}{cc}
\Phi \left( A^{2}\right) & \Phi \left( A\right) \\ 
\Phi \left( A\right) & \Phi \left( A-\log A\right)%
\end{array}%
\right] =\dsum \left[ 
\begin{array}{cc}
\lambda _{j}^{2} & \lambda _{j} \\ 
\lambda _{j} & \lambda _{j}-\log \lambda _{j}%
\end{array}%
\right] \otimes \Phi \left( P_{j}\right) .
\end{equation*}%
For $x>0,$ $x-\log x\geq 1.$ So the Schur product%
\begin{equation*}
\left[ 
\begin{array}{cc}
\lambda _{j}^{2} & \lambda _{j} \\ 
\lambda _{j} & \lambda _{j}-\log \lambda _{j}%
\end{array}%
\right] =\left[ 
\begin{array}{cc}
\lambda _{j}^{2} & \lambda _{j} \\ 
\lambda _{j} & 1%
\end{array}%
\right] \circ \left[ 
\begin{array}{cc}
1 & 1 \\ 
1 & \lambda _{j}-\log \lambda _{j}%
\end{array}%
\right]
\end{equation*}%
is psd. $\blacksquare $\vskip0.1inIn this context, one can easily obtain the
following inequalities%
\begin{equation*}
\left[ 
\begin{array}{cc}
\Phi \left( \log MI-\log A\right) & \Phi \left( \log MA-A\log A\right) \\ 
\Phi \left( \log MA-A\log A\right) & \Phi \left( \log MA^{2}-A^{2}\log
A\right)%
\end{array}%
\right] \geq O
\end{equation*}%
and%
\begin{equation*}
\left[ 
\begin{array}{cc}
\Phi \left( \log A-\log mI\right) & \Phi \left( A\log A-\log mA\right) \\ 
\Phi \left( A\log A-\log mA\right) & \Phi \left( A^{2}\log A-\log
mA^{2}\right)%
\end{array}%
\right] \geq O.
\end{equation*}

\vskip0.1inBounds for eigenvalues in terms of the entries of the matrix have
been studied extensively in literature, $\left[ 13-16\right] .$ A special
case of our Theorem 2.1 gives inequalities related to extreme eigenvalues.

\textbf{Theorem 3.3. }Let $A$ be any \textbf{\ }Hermitian element of $%
\mathbb{M}(n).$ Let $\mu _{\min }$ and $\mu _{\max }$ be the smallest and
largest eigenvalues of $A-\varphi (A)I.$ Then the cubic equation%
\begin{equation}
x^{3}+\frac{\beta _{1}}{\gamma }x^{2}+\frac{\beta _{2}}{\gamma }x+\frac{%
\beta _{3}}{\gamma }=0  \tag{3.6}
\end{equation}%
is positive or negative according as $x=\mu _{\max }$ or $x=\mu _{\min },$
where 
\begin{eqnarray*}
\beta _{1} &=&-\varphi \left( B^{4}\right) \varphi \left( B^{3}\right)
-\left( \varphi \left( B^{2}\right) \right) ^{2}\varphi \left( B^{3}\right)
+\varphi \left( B^{2}\right) \varphi \left( B^{5}\right) \\
\beta _{2} &=&-\varphi \left( B^{3}\right) \varphi \left( B^{5}\right)
+\left( \varphi \left( B^{4}\right) \right) ^{2}+\varphi \left( B^{3}\right)
^{2}\varphi \left( B^{2}\right) -\left( \varphi \left( B^{2}\right) \right)
^{2}\varphi \left( B^{4}\right) \\
\beta _{3} &=&2\varphi \left( B^{2}\right) \varphi \left( B^{3}\right)
\varphi \left( B^{4}\right) -\left( \varphi \left( B^{2}\right) \right)
^{2}\varphi \left( B^{5}\right) -\left( \varphi \left( B^{3}\right) \right)
^{3} \\
\gamma &=&\left( \varphi \left( B^{3}\right) \right) ^{2}-\varphi \left(
B^{2}\right) \varphi \left( B^{4}\right) +\left( \varphi \left( B^{2}\right)
\right) ^{3}.\text{ \ }
\end{eqnarray*}%
\textbf{Proof. }It follows from $\left( 3.3\right) $, that%
\begin{equation}
\left\vert 
\begin{array}{ccc}
-a & \varphi (B^{2}) & \varphi (B^{3})-a\varphi (B^{2}) \\ 
\varphi (B^{2}) & \varphi (B^{3})-a\varphi (B^{2}) & \varphi
(B^{4})-a\varphi (B^{3}) \\ 
\varphi (B^{3})-a\varphi (B^{2}) & \varphi (B^{4})-a\varphi (B^{3}) & 
\varphi (B^{5})-a\varphi (B^{4})%
\end{array}%
\right\vert \geq 0.  \tag{3.7}
\end{equation}%
Also, $\lambda _{\min }I\leq A\leq \lambda _{\max }I.$ So, $a=$ $\lambda
_{\min }-\varphi (A)=$ $\mu _{\min }.$ Expanding the determinant $\left(
3.7\right) $ we see that the expression $\left( 3.6\right) $ is
non-positive. Likewise, the inequality $\left( 3.4\right) $ implies that $%
\left( 3.6\right) $ is non-negative for $b=\mu _{\max }.$ $\ \blacksquare $

\textbf{Example}: Let%
\begin{equation*}
A=\left[ 
\begin{array}{ccc}
3 & -3\sqrt{2} & -9 \\ 
-3\sqrt{2} & -6 & -3\sqrt{2} \\ 
-9 & -3\sqrt{2} & 3%
\end{array}%
\right] .
\end{equation*}%
The estimates of Wolkowicz and Styan $\left[ 16\right] $ gives $\lambda
_{\min }\leq -\sqrt{48}=-6.928$ and $\lambda _{\max }\geq 6.928$ and
estimates of Sharma et al. $\left[ 14\right] $ gives $\lambda _{\min }\leq -%
\sqrt{96}=-9.792$ and $\lambda _{\max }\geq 9.792.$ Our Theorem 3.3 shows
that the $\lambda _{\min }$ is less than equal to smallest root of $%
x^{3}-144x=0.$ So, $\lambda _{\min }\leq -12.$ Likewise, $\lambda _{\max
}\geq 12.$ The eigenvalues of $A$ are $-12,$ $0$ and $12.$

\vskip0.1inWe finally remark that our technique can be extended to study
inequalities involving normal matrices. Choi $\left[ 7,8\right] $ showed
that for any normal element $A$ of $\mathbb{M}\left( n\right) ,$ 
\begin{equation*}
\Phi \left( A\right) \Phi \left( A^{\ast }\right) \leq \Phi \left( A^{\ast
}A\right) ,\text{ \ }\Phi \left( A^{\ast }\right) \Phi \left( A\right) \leq
\Phi \left( A^{\ast }A\right) .
\end{equation*}%
Equivalently%
\begin{equation*}
\left[ 
\begin{array}{cc}
I & \Phi (A) \\ 
\Phi (A^{\ast }) & \Phi (A^{\ast }A)%
\end{array}%
\right] \geq O.
\end{equation*}%
Note that the matrix%
\begin{equation*}
\left[ 
\begin{array}{ccc}
1 & \lambda _{j} & \left\vert \lambda _{j}\right\vert ^{2} \\ 
\overline{\lambda _{j}} & \left\vert \lambda _{j}\right\vert ^{2} & 
\overline{\lambda _{j}}\left\vert \lambda _{j}\right\vert ^{2} \\ 
\left\vert \lambda _{j}\right\vert ^{2} & \lambda _{j}\left\vert \lambda
_{j}\right\vert ^{2} & \left\vert \lambda _{j}\right\vert ^{4}%
\end{array}%
\right]
\end{equation*}%
is psd. Therefore, for any positive unital linear map $\Phi $ and for every
normal matrix $A$ $\in \mathbb{M}\left( n\right) $, we have

\begin{equation}
\left[ 
\begin{array}{ccc}
I & \Phi (A) & \Phi (A^{\ast }A) \\ 
\Phi (A^{\ast }) & \Phi (AA^{\ast }) & \Phi (A^{\ast 2}A) \\ 
\Phi (A^{\ast }A) & \Phi (A^{\ast }A^{2}) & \Phi (A^{\ast 2}A^{2})%
\end{array}%
\right] \geq O.  \tag{3.8}
\end{equation}%
So, for every positive unital linear functional $\varphi $, we have 
\begin{equation}
\varphi \left( \left\vert B\right\vert ^{4}\right) \geq \frac{\left\vert
\varphi \left( B\left\vert B\right\vert ^{2}\right) \right\vert ^{2}}{%
\varphi \left( \left\vert B\right\vert ^{2}\right) }+\left( \varphi \left(
\left\vert B\right\vert ^{2}\right) \right) ^{2},  \tag{3.9}
\end{equation}%
where $B=A-\varphi (A)I$ and $\left\vert B\right\vert ^{2}=B^{\ast }B.$ The
inequality (3.8) is a noncommutative analogue, and $\left( 3.9\right) $ is a
complex analogue of the classical inequality $\left[ 12\right] $%
\begin{equation*}
\mathbb{E}\left[ Y^{4}\right] \geq \frac{\mathbb{E}\left[ Y^{3}\right] ^{2}}{%
\mathbb{E}\left[ Y^{2}\right] }+\mathbb{E}\left[ Y^{2}\right] ^{2},
\end{equation*}%
where $Y=X-\mathbb{E}\left[ X\right] $.

\vskip0.2in\noindent \textbf{Acknowledgements}. The authors are grateful to
Prof. Rajendra Bhatia for the useful discussions and suggestions, and I.S.I.
Delhi for a visit in January 2015 when this work had begun. The support of
the UGC-SAP is also acknowledged.


\begin{thebibliography}{99}
\bibitem{1} R. Bhatia, Matrix Analysis, Springer Verlag New York, 2000.

\bibitem{2} R. Bhatia, C. Davis, A better bound on the variance, Amer. Math.
Monthly 107 (2000) 353-357.

\bibitem{3} R. Bhatia, Positive Definite Matrices, Princeton University Press%
\textit{, }2007.

\bibitem{4} R. Bhatia, R. Sharma, Some inequalities for positive linear
maps, Linear Algebra Appl. 436 (2012) 1562-1571.

\bibitem{5} R. Bhatia, R. Sharma, Positive linear maps and spreads of
matrices, Amer. Math. Monthly 121 (2014) 619-624.

\bibitem{6} R. Bhatia, R. Sharma, Positive linear maps and spreads of
matrices-II, Linear Algebra Appl. 491 (2016) 30-40.

\bibitem{7} M.D. Choi, A Schwarz inequality for positive linear maps on $%
C^{\ast }$- algebras, Illinois J. Math. 18 (1974) 565-574.

\bibitem{8} M.D. Choi, Some assorted inequalities for positive linear maps
on $C^{\ast }$- algebras, J. Operator Theory 4 (1980) 271-285.

\bibitem{9} C. Davis, A Schwarz inequality for convex operator functions,
Proc. Amer. Math. Soc. 8 (1957) 42-44.

\bibitem{10} R.A. Horn, C.R. Johnson, Matrix Analysis, Cambridge University
Press, 2013.

\bibitem{11} R.V. Kadison, A generalized Schwarz inequality and algebraic
invariants for operator algebras, Ann. Math. 56 (1952) 494-503.

\bibitem{12} K. Pearson, Mathematical contributions to the theory of
evolution, XIX: Second supplement to a memoir on skew variation, Philos.
Trans. Roy. Soc. London Ser. 216(A) (1916) 429-457.

\bibitem{13} R. Sharma, Some more inequalities for Arithmetic Mean, Harmonic
Mean and Variance, Journal of Mathematical Inequalities 2(1) (2008\textbf{) }%
109-114.

\bibitem{14} R. Sharma, R. Bhandari, M. Gupta, Inequality related to the
Cauchy-Schwarz inequality, Sankhya 74(A) (2012) 101-111.

\bibitem{15} R. Sharma, R. Kumar, R. Saini, G. kapoor, Bounds on spread of
matrices related to fourth central moment, Bull. Malays. Math. Sci. Soc. DOI
10.1007/540840-015-0267-1.

\bibitem{16} H. Wolkowicz, G.P.H Styan, Bounds for eigenvalues using traces,
Linear Algebra Appl. 29 (1980) 471-506.
\end{thebibliography}
\end{document}